\documentstyle{amsart}
\newtheorem{theorem}{Theorem}[section]
\newtheorem{lemma}[theorem]{Lemma}

\theoremstyle{definition}
\newtheorem{ex}{Example}
\newenvironment{example}{\begin{ex}\rm}{\end{ex}}
\newenvironment{remark}{\m\ni{\it Remark.}\rm}{\m}
\def\c{$C^{\displaystyle *}\!$-}
\def\tensor{\otimes}
\def\CG{C^*\!(G)}
\def\ni{\noindent}
\def\m{\medskip}
\def\Ind{\mathop{\mathrm{Ind}}\nolimits}
\def\by{{\mathrm{\ by\ }}}
\def\id{\mathop{\mathrm{id}}\nolimits}

\def\clsp{\mathop{\overline{\mathrm{sp}}}\nolimits}
\begin{document}

\title{Full crossed products by Hopf \c algebras}
\author{May Nilsen}
\address{Department of Mathematics and Statistics, 
      University of Nebraska-Lincoln,
      Lincoln NE 68588-0323}
\email{mnilsen@@math.unl.edu}


\begin{abstract}
We show that when a co-involutive Hopf
 \c algebra $S$ coacts via $\delta$ on a \c algebra $A$,
 there exists a full crossed product
 $A\times_\delta S$,
with universal properties analogous to those
of full crossed products by locally compact groups.
The dual Hopf \c algebra is then defined by $\hat S:={\Bbb C}\times_{\id} S$.
\end{abstract}

\maketitle
 
\section*{introduction.}\label{sec:intro}

Due to differences in notation and construction, it is often considered that
 coactions of locally compact groups on \c algebras
 lead to a  type of crossed product different
to that of  actions of locally compact groups.
But the definitions of the two 
types of full crossed product by their universal properties given by
Raeburn \cite{iain:tak,iain:coact} 
are quite similar.
With a slight change of notation, actions and coactions 
can look very similar indeed. 
By taking the point of view that these   crossed products
are essentially of the same type,
we were able to give a short proof of the Imai-Takai duality theorem \cite{it},
the key ideas of which also worked for coactions, 
thus giving a short proof of Katayama's duality theorem \cite{kat} as well 
\cite{nil:dual}.

The natural umbrella under which both actions and coactions fit
is that of co-involutive Hopf \c algebras.
Both the full group \c algebra  $\CG$ and the
algebra of continuous functions on the group disappearing at infinity $C_0(G)$ are examples
of co-involutive Hopf \c algebras. 
By using universal properties to define full crossed products of
\c algebras by coactions of   Hopf \c algebras,
one should be able to obtain full crossed products by locally compact groups by
actions and coactions as two special cases.
 Ng suggested  that it should be possible to 
do this using  
 Raeburn's method
 \cite[Prop 2.13]{ng:cc}.

That is what we accomplish in this paper. 
We will define a class of Hopf \c algebras, called co-involutive,
for which full crossed products can be defined.
At each stage of the construction we look at the examples 
$\CG$ and $C_0(G)$,
 so that it is evident that our construction does give the 
usual full crossed products in those cases.
We also give examples which do not arise from locally compact groups, 
namely those arising from certain amenable multiplicative
unitaries of Baaj and Skandalis \cite{bsk},
so that
it is clear that the full crossed products constructed here
are a genuine generalization of those by  locally compact groups.
The fundamental difference between our crossed products and 
those of  Baaj and Skandalis 
is that ours generalize {\sl full} crossed products,
whereas theirs generalize {\sl reduced} 
crossed products \cite[p11]{baaj}.

We simply define the dual of a co-involutive Hopf \c algebra $S$
to be the trivial crossed product of the complex numbers by $S$.
Using its universal properties we can  readily write down its
comultiplication,
making the dual object a Hopf \c algebra as well.
An integral part of our method is the
construction of a distinguished pair $(\mu_S, V_S)$,
where $\mu_S$ is a representation of $S$, and $V_S$ is a corepresentation of $S$,
on a Hilbert space $H_S$.
They are the analogue of the left regular representation of $\CG$ and
the representation of $C_0(G)$ by multiplication operators on $L^2(G)$.

In \S \ref{sec:group},
we give some notation and summarize   preliminaries concerning locally compact groups.
In \S \ref{sec:hopf}
we define Hopf \c algebras and show 
  that $\CG$ and $C_0(G)$ are Hopf \c algebras;
we also provide other examples.
In \S \ref{subsec:coact},
we define and give examples of coactions of
Hopf \c algebras.
We define co-involutive Hopf \c algebras in \S \ref{sec:dual.hopf},
which are the class of Hopf \c algebras to 
which we can associate  dual Hopf \c algebras.
In this section   a *-algebra $A(S)$ is defined.
In \S \ref{sec:crosprod}, we go on to define   full
 crossed products
by co-involutive Hopf \c algebras.
The dual Hopf \c algebra $\hat S$ of a given 
co-involutive Hopf \c algebra $S$
is defined in
 \S \ref{subsec:can.unit}
and we show that $A(S)$ is dense in  $\hat S$.
Finally, we close with a discussion of the problems surrounding
whether or not the double dual $\widehat{\hat S}$
is isomorphic to $S$.
These, and other issues, will be dealt with in 
\cite{nil:hopfresind}.

\m

I would like to thank   
Etienne Blanchard for his generous assistance 
concerning understanding various parts of \cite{bsk}.
A special thanks goes to my colleague David Pitts at the 
University of Nebraska-Lincoln,
for his substantial help, without which,
 this paper would have been much the poorer.


\section{Preliminaries.}\label{sec:group}

All representations of a \c algebra $A$ are involutive and nondegenerate,
 and so extend uniquely to   unital *-representations of 
the multiplier algebra $M(A)$  \cite[Lem 1.1]{lprs}.
All tensor products  $\tensor$ will be minimal
and  
all identity maps are denoted by $\id$.

For each $f$ in the dual space $B^*$  of a \c algebra $B$, 
there is a well-defined map, called a  slice map,
 $\id\tensor f\colon A\tensor B\to A$ with the property that
$\id\tensor f(a\tensor b)=af(b)$.
Using the Cohen factorization theorem
one can show that slice maps extend to multiplier algebras
 \cite[p115]{man};
the extension is also denoted by $\id\tensor f$.
  Tomiyama showed that slice maps $\id\tensor f$ separate points of $A\tensor B$
\cite[Th 1]{tom:ten}.
When  extended to  
$M(A\tensor  B)$, they still separate points 
 because 
the multiplier algebra $M(A\tensor B)$   is the largest \c algebra containing 
$A\tensor  B$
as an essential ideal \cite[p82]{mur}.

Let $C$ be a \c algebra. 
   Given a   nondegenerate homomorphism $\phi\colon A\to M(C)$,
 there is nondegenerate homomorphism
$\phi\tensor\id\colon A\tensor B\to M(C\tensor B)$ 
\cite[Lem 1.1]{iain:coact}. 
The flip isomorphism
$\Sigma\colon A\tensor B\to B\tensor A$ satisfies 
$\Sigma(a\tensor b)=b\tensor a$.


\m

Let $G$ be a locally compact group,
$\CG$ be the full group \c algebra with canonical embedding 
$i_G\colon G\to \CG$,
and $C_0(G)$ be the continuous functions on $G$ disappearing at infinity.
The Fourier-Stieltjes algebra $B(G)$ of $G$ is the
set of all finite, complex-linear combinations of 
continuous, positive definite
functions on $G$, with pointwise operations
\cite[p21]{marty:w*} \cite[32.10]{hr:II}.
Eymard showed that $B(G)$ can be identified with $\CG^*$ the space of 
linear functionals on $\CG$,
where a   functional $f\in \CG^*$ corresponds to the
function on $G$ defined by     $s\mapsto f(i_G(s))$   \cite[p192]{eymard}.
This means that $\CG^*$ can be considered to be contained in
$C_b(G)$.
Eymard also showed that 
$B(G)$ is a Banach *-algebra
with the norm inherited from $\CG^*$.
 
The Fourier  algebra $A(G)$ is the *-ideal in $B(G)$,
which is the norm closure of the set of functions of compact support   
\cite[p21]{marty:w*}.
The spectrum of $A(G)$ is $G$  \cite[3.34]{eymard},
hence $A(G)$ is a  dense *-subalgebra of $C_0(G)$
in the supremum norm.
 
 The map $w_G :=s\mapsto i_G(s)$ is a continuous bounded function on $G$
into $M(\CG)$,
and since $C_b(G,M(\CG))$ is contained in
$M(C_0(G)\tensor\CG)$ \cite[p751]{lprs},
$w_G$ can be considered to be a unitary element of 
$M(C_0(G)\tensor\CG)$.
The unitary $v_G:=\Sigma(w_G)$ is an element of
  $M(\CG\tensor C_0(G))$.
Note that if $f\in A(G)$, then 
\begin{equation}\label{eq:rock}
\id\tensor f(w_G)
=\id\tensor f(s\mapsto i_G(s))=s\mapsto f(i_G(s))=f\in C_0(G).
\end{equation}
\ni If $f\in \CG^*$, then $\id\tensor f(w_G)=f\in M(C_0(G))$.
For each $t\in G$, 
there is *-homomorphism $\varepsilon_t\colon C_0(G)\to{\Bbb C}$,
evaluation at $t$. It is an element of
$C_0(G)^*$ and
\vspace*{-.3cm}
\begin{equation}\label{eq:rock2}
\id\tensor \varepsilon_t(v_G)
=\id\tensor \varepsilon_t(s\mapsto i_G(s))=i_G(t)\in M(\CG).
\end{equation}
\ni Also, for $z\in L^1(G)\subseteq M1(G)=C_0(G)^*$ 
\cite[7.1.2]{pedbook},
$\id\tensor z(v_G)=z\in  \CG$.

\m

Let $\alpha$ be a strongly continuous action of
a locally compact group $G$    on a 
\c algebra $A$. 
Because each map $s\mapsto\alpha_s(a)$ is   continuous,
they are in $C_b(G, A)$, and making use of the embedding of 
$C_b(G,A)$ in $M(A\tensor C_0(G))$,
we can consider $\alpha$ to be a *-homomorphism from
$A$ to $M(A\tensor C_0(G))$.
The fact that the original $\alpha$ was a group homomorphism,
gives   an identity which the new $\alpha$ satisfies:
\[(\alpha\tensor\id)\circ\alpha
=(\id\tensor\alpha_G)\circ\alpha,\]
\ni as maps on  $M(A\tensor C_0(G)\tensor C_0(G))$, and
where 
 $\alpha_G\colon C_0(G)\to M(C_0(G)\tensor C_0(G))$ is 
the homomorphism  given  by
 \begin{equation}\label{eq:alpha.defn}
\alpha_G(f)(s,t)= f(st).
\end{equation}
\ni This map    satisfies the   comultiplication identity:
 \begin{equation}\label{eq:alpha.comult.identity}
(\alpha_G\tensor\id)\circ\alpha_G
=(\id\tensor\alpha_G)\circ\alpha_G,
\end{equation}
\ni as maps on  $M(C_0(G)\tensor C_0(G)\tensor C_0(G))$
\cite[Rem 2.2(1)]{iain:coact}.
Having expressed   actions of $G$ on  $A$ entirely in terms of
$C_0(G)$,
to define coactions of $G$, we simply replace $C_0(G)$ by $\CG$. 
An   injective  *-homomorphism  
$\delta\colon A\to M(A\tensor \CG)$
  is a   {\it coaction}  of $G$ on $A$ if 
it satisfies the   coaction identity:
$(\delta\tensor\id)\circ\delta
=(\id\tensor\delta_G)\circ\delta,
$
\ni where 
 $\delta_G\colon \CG\to M(\CG\tensor \CG)$ is the 
  injective  *-homomorphism satisfying
\begin{equation}\label{eq:delta.defn}
\delta_G(i_G(s))=i_G(s)\tensor i_G(s).
\end{equation}
\ni 
The map $\delta_G$   satisfies a  comultiplication identity
 like Equation (\ref{eq:alpha.comult.identity})
\cite[p628]{iain:coact}.
A system where $G$ coacts on a \c algebra $A$ via $\delta$ 
is called a   {\it cosystem}  \cite{lprs,nil:coact,qui:frc}.

In the special case where 
 $G$ is abelian,
Pontryagin duality says that the   characters form 
a locally compact abelian group $\widehat G$ 
\cite[Th 1.7.2]{rudin:gp}.
The  Fourier transform
gives an isomorphism between $\CG$ and $C_0(\widehat G)$,
which 
 carries $\delta_G$
to $\alpha_{\widehat G}$.
Thus even when $G$ is
non-abelian, 
  $C_0(G)$ can be thought of as  the dual of $\CG$.


\section{Hopf \c algebras.}\label{sec:hopf}

A  {\it Hopf \c algebra} is a \c algebra $S$ together with
a 
 nondegenerate injective *-homomorphism
$\delta_S\colon S\to M( S\tensor S)$
  such that 
\begin{equation}\label{coact.identity}
(\delta_S\tensor\id)\circ\delta_S
=(\id\tensor\delta_S)\circ\delta_S,
\end{equation}
\ni as maps of $S$ into  $M(S\tensor S\tensor S)$.
The homomorphism $\delta_S$ is called the {\it comultiplication} on $S$
and  Equation (\ref{coact.identity}) is called the {\it comultiplication identity}.
A  Hopf \c algebra  is 
  {\it right-simplifiable} if 
\begin{equation}\label{eq:nondeg}
\delta_S(S)(1 \tensor S)=S\tensor S,
\end{equation}
\ni    
and {\it left-simplifiable} if $\delta_S(S)(S \tensor 1)=S\tensor S$.
A  Hopf \c algebra  is 
{\it bi-simplifiable} if it is both left and right-simplifiable
\cite[Defn 0.1]{bsk}.
%


  A {\it corepresentation} of $S$ on a Hilbert space $H$ 
is a unitary $V\in M(K(H)\tensor S)$ such that
\vspace{-.3cm}
\begin{equation}\label{corep.id}
\id \tensor \delta_S(V)=V_{12}V_{13} 
        \in M(K(H)\tensor S \tensor S),
\end{equation}
\ni  and where $V_{12}$ is the element of 
$M(K(H)\tensor S\tensor S)$ that is $V$ acting on the first
and second factors (that is $V\tensor 1$),
and where $V_{13}$   is $V$ acting on the first
and third factors \cite[Defn 0.3]{bsk}.
Such a unitary $V$   also satisfies
\[
\delta_S\tensor \id (\Sigma V)       
  =(\Sigma V)_{13}(\Sigma V)_{23}     
      \in M(S \tensor S\tensor K(H)).
\]
%


\m\ni{\bf Example 1(a).}
Let $A$ be any \c algebra.
It can be given the trivial right-simplifiable comultiplication
$\delta\colon A \to M(A\tensor A)$ satisfying
$\delta(a)=a \tensor 1$.
The unitary  $1\tensor 1\in M(K(H)\tensor A)$  is the only 
corepresentation because
$\id\tensor\delta(U)=U_{12}$ 
  must equal $U_{12}U_{13}$.


\m\ni{\bf Example 2(a).} 
The \c algebra $C_0(G)$ is   a Hopf \c algebra
with right-simplifiable comultiplication  
 $\alpha_G\colon C_0(G)\to M(C_0(G)\tensor C_0(G))$ 
(Equation (\ref{eq:alpha.defn}))  
\cite[Rem  2.2(2)]{iain:coact}.
Let   $\lambda$ be the left regular representation of $\CG$ on $L^2(G)$. 
Since $\lambda\colon \CG\to M(K(L^2(G)))$,
 \cite[Lem 1.1]{iain:coact}
implies that the unitary
$\lambda\tensor\id(v_G)$
is    in $M(K(L^2(G))\tensor  C_0(G))$.
The following shows that
  $\lambda\tensor\id(v_G)$  is a  corepresentation of $C_0(G)$:
\begin{eqnarray*}
\id\tensor\alpha_G(\lambda\tensor\id(v_G))
&=&\lambda\tensor\id\tensor \id(\id\tensor\alpha_G(v_G))\\
&=&\lambda \tensor\id\tensor\id \lambda((s,t)\mapsto i_G(st))\\
&=&\lambda \tensor\id\tensor\id \lambda((s,t)\mapsto i_G(s)i_G(t))\\
&=&(\lambda\tensor\id(v_G))_{12}(\lambda\tensor\id(v_G))_{13},
\end{eqnarray*}
\ni in $M(K(L^2(G))\tensor C_0(G)\tensor C_0(G))$.
\ni Any *-representation $ U\colon \CG\to B(H)$ gives rise to a 
 corepresentation of $C_0(G)$, namely
       $U\tensor\id(v_G)\in M(K(H)\tensor C_0(G))$.

 Conversely, a corepresentation $V$ of $C_0(G)$ is a unitary element
of \newline $M(K(H)\tensor C_0(G))\subseteq C_b(G,B(H))$,
and so $V$ corresponds to a  continuous bounded function $\tilde V$
from $G$ into the unitary group of $B(H)$.
We know that $\id\tensor \alpha_G(V)=V_{12}V_{13}$.
The left hand side is the function
on $G\times G$ that takes
$(s,t)$ to $\tilde V(st)$,
and the right hand side is the function
  that takes
$(s,t)$ to $\tilde V(s)\tilde V(t)$.
This means that $\tilde V$  is a unitary group representation
and so integrates up to a *-representation $\bar V$ of $\CG$.
To know  that all  corepresentations of $C_0(G)$ are of the form
$\bar V\tensor\id(v_G)$ for some *-representation $\bar V$ of
 $\CG$,
we need to show that $\bar V\tensor\id(v_G)=V$.
The point evaluation functionals $\varepsilon_s$
on $C_0(G)$ are sufficient to separate points.
So functionals of the form $\id\tensor\varepsilon_s$ separate points
of $M(K(H)\tensor C_0(G))$ \cite[Th 1]{tom:ten} and
  the following calculation suffices:
\[\id\tensor\varepsilon_s(\bar V\tensor\id(v_G))
= \bar V(\id\tensor\varepsilon_s(v_G))
= \tilde V(s)= \id\tensor\varepsilon_s(V).\]


\m\ni{\bf Example 3(a).} 
In \cite[Ex 2.3(1)]{iain:coact}
Raeburn notes that
  $\delta_G$ 
(Equation (\ref{eq:delta.defn}))
 is a right-simplifiable comultiplication on $\CG$.
Let $M$ be the representation of $C_0(G)$ on $L^2(G)$
as multiplication operators.
 Then
$M\tensor\id(w_G)$
 is a  \hbox{corepresentation of $\CG$:}
\begin{eqnarray*}
\id\tensor\delta_G(M\tensor\id(w_G))
&=&M\tensor\id\tensor \id(\id\tensor\delta_G(w_G))\\
&=&M\tensor\id\tensor \id(s\mapsto i_G(s)\tensor i_G(s))\\
&=&M\tensor\id\tensor \id(s\mapsto i_G(s)\tensor 1
       .\,s\mapsto 1\tensor i_G(s))\\
&=&M\tensor\id(w_G)\tensor 1
       .\id\tensor\Sigma(M\tensor\id(w_G)\tensor 1))\\
&=&(M\tensor\id(w_G))_{12}(M\tensor\id(w_G))_{13}.
\end{eqnarray*}
\ni Any *-representation 
  $\mu\colon C_0(G)\to B(H)$ gives rise to a  
        corepresentation of $\CG$, namely 
         $\mu\tensor\id(w_G)\in M(K(H)\tensor\CG)$. 

Now for the converse.
Given a corepresentation 
 $U$ of $\CG$ in $M(K(H)\tensor\CG)$,
we can define a map $\bar\mu$ on  $\CG^*$  to $B(H)$ by
$\bar\mu(f):=\id\tensor f(U).$
\ni We denote its restriction to the
 Fourier algebra $A(G)$ by $\mu$.
We first show that $\mu$ is multiplicative:
\begin{eqnarray*}
\mu(fg)
&=&\id \tensor (fg)(U) 
 = \id\tensor ((f\tensor g)\circ\delta_G)(U)\\
&=&\id\tensor f\tensor g(\id\tensor\delta_G(U))
=\id\tensor f\tensor g(U_{12}U_{13})\\
&=&\id\tensor f(U)\id\tensor g(U)
=\mu(f)\mu(g).
\end{eqnarray*}
\ni  
An argument   described in \cite[p130]{nt}
 shows that $\mu$ is involutive. 
Thus $\mu$ is a *-homomorphism from the 
abelian Banach *-algebra $A(G)$ into $B(H)$.
The fact that it   is *-preserving implies that it is norm decreasing on $A(G)$
\cite[Th 2.1.7]{mur}.
Thus there exists a unique *-representation,
also denoted by $\mu$, on the   enveloping \c algebra 
\cite[2.7.4]{dixbook}.
Since the spectrum of $A(G)$ is $G$  \cite[3.34]{eymard}, 
the  enveloping \c algebra is $C_0(G)$.

To see that  every corepresentation of $\CG$
is of the form
$\mu\tensor\id(w_G)$,
we   show that 
  $\mu\tensor\id(w_G)=U$. 
\ni By definition 
$\bar\mu(f)=\id\tensor f(U)$ for all $f\in \CG^*$.
Using Equation (\ref{eq:rock}),
we have
$\mu(\id\tensor f(w_G))=\id\tensor f(U)$,
which means
$\id\tensor f(\mu\tensor\id(w_G))=\id\tensor f(U)$.
Since  
functionals of the form $\id\tensor f$   
 separate the points of $M(K(H)\tensor\CG)$ \cite[Th 1]{tom:ten},
 we deduce that $\mu\tensor\id(w_G)=U$
(cf. \cite[Lem 1.2]{qr:ind}).


\m\ni{\bf Example 4(a).}
Let   $\alpha$ be an action of $G$   on $A$.
There is a dual coaction $\hat\alpha$ on the full crossed product
$A\times_\alpha G$.
Specifically,
$\hat\alpha\colon A\times_\alpha G
   \to M((A\times_\alpha G)\tensor \CG)$ such that
$\hat\alpha(k_A(a))=k_A(a)\tensor 1$ and 
$\hat\alpha(k_G(s))=k_G(s)\tensor i_G(s)$
 \cite[Ex 2.3(1)]{iain:coact}.
Composing this map with 
$\id\tensor k_G$, where  $k_G$ is the canonical homomorphism 
carrying $\CG$ into $A\times_\alpha G$,
gives a map 
$\delta_{A\times G}\colon
A\times_\alpha G \to M((A\times_\alpha G)\tensor (A\times_\alpha G))$
such that
\[\delta_{A\times G}(k_A(a))=k_A(a)\tensor 1 \ \and \ 
 \delta_{A\times G}(k_G(s))=k_G \tensor k_G(\delta_G(s)).\]
\ni This is
 a right-simplifiable comultiplication on $A\times_\alpha G$,
so $A\times_\alpha G$ is a Hopf \c algebra.

Let $W$ be a corepresentation of $\CG$.
The following shows that $\id\tensor k_G(W)$
is a 
  corepresentation of $A\times_\alpha G$:
\begin{eqnarray*}\id\tensor\delta_{A\times G}(\id\tensor k_G(W))
&=&\id\tensor k_G\tensor k_G(\id\tensor \delta_G(W))\\
&=&\id\tensor k_G\tensor k_G(W_{12}W_{13})\\
&=&\id\tensor k_G\tensor k_G(W_{12})\id\tensor k_G\tensor k_G(W_{13})\\
&=&\id\tensor k_G(W)_{12}\id\tensor k_G(W)_{13}.
\end{eqnarray*}
As yet we have not been able to show that all corepresentations of 
$A\times_\delta G$ are of this form,
although we expect this to be the case.

In the case where $A$ is abelian,
the crossed product is   a transformation group \c algebra.
These \c algebras are a type of groupoid \c algebra.
So we have shown   that at least a subclass
of groupoid \c algebras are Hopf \c algebras.

In the case where $G$ is acting by right translation $\sigma$ on $C_0(G)$,
the crossed product $C_0(G)\times_\sigma G$
is isomorphic to $K(L^2(G))$ \cite{rie:compacts},
so we have that $K(L^2(G))$ is a Hopf \c algebra.
When $G$ is a finite group of order $n$, this leads to 
identifying the $n\times n$ matrices $M_n(\Bbb{C})$ as 
a Hopf \c algebra.


\m\ni{\bf Example 5(a).}
Let $\delta$ be a coaction of a locally compact group $G$   a
\c algebra $A$.
There is a dual  action $\hat\delta$ on the full crossed product
$A\times_\delta G$ \cite[Cor 2.14]{iain:coact}.
We can compose this with $\id\tensor j_{C_0(G)}$
to  obtain a map $\delta_{A\times G}$ from
$A\times_\delta G$ into  $M((A\times_\delta G)\tensor (A\times_\delta G))$
making $A\times_\delta G$ into  a Hopf \c algebra.
As in Example 4(a), if
 $W$ is a corepresentation of $C_0(G)$,
then $\id\tensor j_{C_0(G)}(W)$
is a   corepresentation of $A\times_\delta G$.


\m\ni{\bf Example 6(a).} 
Let $V\in B(H\tensor H)$ be a regular multiplicative unitary 
on a Hilbert space $H$, as defined by
Baaj and Skandalis \cite[Defn 1.1, 3.3]{bsk}.
This implies, among other things,
 that $V$ is a unitary operator which satisfies the 
pentagonal identity:
\[V_{12}V_{13}V_{23}=V_{23}V_{12}.\]
\ni In \cite[Cor A.6]{bsk} they associate to each such $V$ a pair
of Hopf \c algebras, $S_p$ and $\hat S_p$.
They also define a pair
of Hopf \c algebras, $S_V:=\overline{A(V)}$ and $\hat S_V:=\overline{\hat A(V)}$,
where 
\[A(V):={\rm sp}\{f\tensor\id(V)|f\in B(H)_*\}
\ \and \ 
 \hat A(V):={\rm sp}\{\id\tensor f(V)|f\in B(H)_*\},\]

\ni \cite[Defn 1.3, 1.5]{bsk}, 
which should be considered as the reduced version of the theory.

Baaj and Skandalis show  \cite[\S 4]{bsk},
that the compact quantum groups of Woronowicz \cite{wor:cmp} give rise to
compact regular multiplicative unitaries. 
This    provides
 examples of Hopf \c algebras which do not arise 
from locally compact groups.


\addtocounter{ex}{6}

\m

The following examples are not actually examples --
 they are examples   which one might think could be 
Hopf \c algebras, but aren't.

\begin{example}\label{ex:twist}
Let $(\id, u)$ be a twisted action of a locally compact group
$G$ 
on the complex numbers $\Bbb{C}$ \cite[Defn 2.1]{PRI}.
We recall that the twist $u$ is a Borel map from $G\times G$ to $\Bbb{T}$
satisfying $i_G(s)i_G(t)=u(s,t)i_G(st)$.
 The crossed product ${\Bbb C}\times_{\id,u}G$
for this twisted dynamical system is called a 
{\it twisted group \c algebra},
sometimes denoted by $C^{\displaystyle *}\!(G,u)$.
We would like to define a 
 comultiplication 
 $\delta_G
   \colon C^{\displaystyle *}\!(G,u)\to 
   M(C^{\displaystyle *}\!(G,u)\tensor 
     C^{\displaystyle *}\!(G,u))$ 
as we did in Example 2(a),
but   $\delta_G$ is {\sl not} a homomorphism:
\begin{eqnarray*}
\delta_G(i_G(s)i_G(t))
&=&\delta_G(u(s,t)i_G(st))\\
&=& u(s,t)\delta_G(i_G(st))
= u(s,t)(i_G(st)\tensor i_G(st)),\ \and
\end{eqnarray*}
\vspace{-0.5cm}
\begin{eqnarray*}
\delta_G(i_G(s))\delta_G(i_G(t))
&=&\delta_G(i_G(s))\delta_G(i_G(t))\\
&=&[i_G(s)\tensor i_G(s)]\,[i_G(t)\tensor i_G(t)]\\
&=&i_G(s)i_G(t)\tensor i_G(s)i_G(t) 
= u(s,t)i_G(st)\tensor u(s,t)i_G(st). 
\end{eqnarray*}
\ni So $C^{\displaystyle *}\!(G,u)$ cannot be made into a Hopf \c algebra
in this manner. 
\end{example}

\begin{example}\label{ex:semi}
Let $\Gamma$ be a totally ordered discrete abelian group
with positive cone $\Gamma^+$.
The semigroup \c algebra  
 $C^{\displaystyle *}\!(\Gamma^+)$ 
has been defined \cite[\S 1]{alnr}
and is the universal \c algebra whose representations
correspond to the isometric representations of the semigroup
$\Gamma^+$.
Routine calculations show that we can define a 
  comultiplication  $\delta_{\Gamma^+}$
as in Example 3(a),
so  $C^{\displaystyle *}\!(\Gamma^+)$  
is a Hopf \c algebra.
We can again define 
\newline
$w_{\Gamma^+}\in M(C_0(\Gamma^+)\tensor C^{\displaystyle *}\!(\Gamma^+))$,
but since each $i_{\Gamma^+}\!(s)$ is a not necessarily unitary isometry,
$w_{\Gamma^+}$ may not be unitary,
and so elements of the form $\mu\tensor\id(w_{\Gamma^+})$
are not 
 corepresentations of $C^{\displaystyle *}\!(\Gamma^+)$
 in the sense we are using in this paper.

The function algebra $C_0(\Gamma^+)$ has comultiplication 
$\alpha_{\Gamma^+}$
as in Example 2(a),
but again it does not necessarily have any unitary corepresentations.
\end{example}


\begin{example}\label{ex:groupoid}
Let ${\cal G}$ be a groupoid.
In Example 4(a) we showed how   groupoid \c algebras 
  arising as transformation group \c algebras 
are Hopf \c algebras.
But in general 
there does not seem to be an obvious way to 
make a groupoid \c algebra $C^{\displaystyle *}\!({\cal G})$
into a Hopf \c algebra.
 Working with von Neumann algebras, 
Yamanouchi  \cite[p12]{yam:groupoid}  shows that 
  $L^\infty({\cal G})$ 
admits a comultiplication $\alpha_{\cal G}$ satisfying 
\[\alpha_{\cal G}(f)(s,t)=\left\{ \begin{array}{ll}
                       f(st),& \mbox{if $(s,t)$ is a composable pair,}\\
                       0,& \mbox{otherwise.}\\
                       \end{array}     
                     \right. \]
\ni This approach does not work for $C_0({\cal G})$,
as $\alpha_{\cal G}(f)$ need not be continuous on ${\cal G}\times {\cal G}$.
\end{example}

\begin{remark}
Another concept frequently seen in the literature is that of a co-unit
\cite[Defn 1.5(b)]{ng:cc}.
We have no use for a co-unit in this paper.
Actually,  there is a problem with the name ``co-unit".
It's logical to think that a co-unit is something that would be used
to define a unit in the dual object.
But consider the fundamental example: $\CG$  and its dual $C_0(G)$.
The trivial representation of $G$ integrates up to a representation of $\CG$,
denoted by $1$, and
$(\id\tensor 1)\circ\delta_G=\id$
\ni \cite[Lem 1.3]{iain:coact}, 
which means $\CG$ is  co-unital.
However, unless the group is compact, $C_0(G)$ is not unital.
So, the existence of such a co-unit does not imply that the dual is unital.
\end{remark}



\section{Coactions of Hopf \c algebras.}\label{subsec:coact}

 A {\it  coaction} of  a Hopf \c algebra $S$ on a \c algebra $A$ 
is a nondegenerate injective  *-homomorphism 
$\delta\colon A\to M(A\tensor S)$
   such that
$(\delta\tensor \id)\circ\delta=(\id\tensor\delta_S)\circ\delta$,
\ni as maps of $A$ into  $M(A\tensor S\tensor S)$.
If 
 $\delta(A)(1\tensor S)= A\tensor S$, then $\delta$ is called a 
{\it nondegenerate} coaction 
(cf. \cite[Defn 0.2]{bsk}, \cite[Defn 1.2]{lprs}).

Let  $B$ be a \c algebra.
A {\it covariant homomorphism} of $(A,S,\delta)$ into 
$M(B)$ is a pair $(\phi, v)$, 
where $\phi\colon A\to M(B)$ is a nondegenerate  *-homomorphism and 
 $v\in M(B\tensor S)$ is a unitary  such that, 
\begin{itemize}
\item[(a)] $\id \tensor \delta_S(v)=v_{12}v_{13} 
   \in M(B\tensor S \tensor S)$, and 
\item[(b)] $\phi\tensor \id(\delta(a))
     =v(\phi(a)\tensor 1)v^*\in M(B\tensor S)$, 
                       for all  $a\in A$. 
\end{itemize}
\ni A {\it covariant representation} of  $(A, S, \delta)$ 
is a   covariant homomorphism of $(A,S,\delta)$ into $M(K(H))\cong B(H)$.
Specifically, it is a pair $(\pi,V)$, 
where $\pi\colon A\to B(H)$ is a *-representation of $A$ 
and $V\in  M(K(H)\tensor S)$ is a  corepresentation of $S$,
such that
\[\pi\tensor \id(\delta(a)) =V(\pi(a)\tensor 1)V^* \in M(K(H)\tensor S),\]
\ni  \cite[Defn 0.2, 0.3]{bsk}.
Note that this implies that 
$\ker((\pi\tensor \id)\circ\delta)=\ker(\pi)$.
 

\begin{example}\label{ex:hopf.C^*-alg=hopf.system}
 There is always the trivial coaction of a Hopf \c algebra $S$
on  a \c algebra $A$,
that is,  $\id\colon A\to  M(A\tensor S)$
satisfying $\id(a)=a\tensor 1$.
Also, a right-simplifiable 
comultiplication $\delta_S\colon S\to M(S\tensor S)$ 
is a coaction of a Hopf \c algebra $S$ on itself. 
\end{example}

\m\ni{\bf Example 4(b).}
Let $(A, G, \alpha)$ be a dynamical system,
so that, as described in \S 1,
we have the map 
 $\alpha\colon A\to M(A\tensor C_0(G))$.
This map
  is a coaction of the Hopf \c algebra $C_0(G)$ on $A$.
A pair $(\pi,U)$    which is a covariant representation of
 the dynamical system $(A, G, \alpha)$,
that is, $\pi\colon A\to B(H)$ and $U\colon \CG\to B(H)$ 
are representations such  that
\[\pi\tensor\id(\alpha(a))
=U\tensor\id(v_G)(\pi(a)\tensor 1)U\tensor\id(v_G^*)
        \ \in\  M(K(H)\tensor C_0(G)),\]
\ni gives rise to a covariant pair 
$(\pi, U\tensor\id(v_G))$ of 
   $(A, C_0(G), \alpha)$, because
  the covariance condition is precisely the same,
and, as we showed in Example 2(a),  
   $U\tensor\id(v_G)$
 is a corepresentation of $C_0(G)$.


\m\ni{\bf Example 5(b).}
Let $(A, G, \delta)$ be a cosystem.
Then $\delta$ is a coaction of the Hopf \c algebra $\CG$ on $A$.
A pair $(\pi,\mu)$ which is a covariant representation of 
the cosystem $(A, G, \delta)$,
that is, $\pi\colon A\to B(H)$ and $\mu\colon C_0(G)\to B(H)$ 
are representations such  that
\[\pi\tensor\id(\delta(a))
=\mu\tensor\id(w_G)(\pi(a)\tensor 1)\mu\tensor\id(w_G^*)
        \ \in\  M(K(H)\tensor\CG),\]
\ni gives rise to a covariant pair 
$(\pi, \mu\tensor\id(w_G))$ of 
 $(A, \CG, \delta)$.


\m\ni{\bf Example 6(b).}
Baaj and Skandalis point out that 
if $(\pi,V)$ is a covariant representation for any
coaction of a Hopf \c algebra of a \c algebra $A$,
then $\id\tensor\pi(V)$ is a multiplicative unitary
\cite[Ex 1.2(5)]{bsk}.
This is not difficult to see because the covariance condition implies that
\[\id\tensor \pi\tensor \id(\id \tensor\delta(V))(1\tensor V) 
=(1\tensor V)(\id\tensor \pi(V)\tensor 1),\]
\ni so using Equation (\ref{corep.id}) gives
$\id\tensor \pi\tensor \id(V_{12}V_{13})V_{23} 
=V_{23}(\id\tensor \pi\tensor \id(V_{12})).$
\ni By applying $(\id\tensor \id\tensor \pi)$ to both sides
we obtain
\[(\id\tensor \pi(V))_{12}\, (\id\tensor \pi(V))_{13}\,(\id\tensor \pi(V))_{23} 
=(\id\tensor \pi(V))_{23}\,(\id\tensor \pi(V))_{12}.\]
%


\section{Co-involutive Hopf \c algebras.}\label{sec:dual.hopf}
 
Let $\delta$ be a coaction of a Hopf \c algebra $S$ on a \c algebra $A$.
 Let $(\pi_1, W_1)$  and $(\pi_2, W_2)$
be   covariant representations of $(A, S,\delta)$ on $H_1$ and $H_2$
respectively.
The direct sum $\pi_1\oplus\pi_2$ is a representation
of $A$ on $H_1\oplus H_2$ and
the unitary $W_1\oplus W_2$ can be
considered an element of
$M(K(H_1\oplus H_2)\tensor S)$.
Furthermore,  $W_1\oplus W_2$ is a corepresentation of $S$ 
and   $(\pi_1\oplus\pi_2,W_1\oplus W_2)$
is a covariant representation of $(A, S,\delta)$ on $H_1\oplus H_2$.
Infinite direct sums are handled similarly.
All these follow from Raeburn's arguments for the group case
\cite[p634]{iain:coact}. 
Now let $W_1$ and $W_2$ be corepresentations of $S$ on $H$.
Then $W_1$ is {\it unitarily equivalent} to $W_2$
if there exists a unitary $U$ in $B(H)$ such that
$W_1=(U\tensor 1)W_2(U^*\tensor 1)$
\hbox{\cite[Eq5.2]{wor:twi}.}

\begin{remark}
Example \ref{ex:semi} suggests that it is possible to have
Hopf \c algebras which do not have any corepresentations,
and consequently no covariant representations.
From here on we will only work with Hopf \c algebras that admit at least
one covariant representation. 
\end{remark}

Let $S$ be a   Hopf \c algebra.
Let $\Gamma$ be the smallest set of  covariant representations  
of   $(S, S,\delta_S)$
such that, for every  covariant representation $(\nu,W)$ of 
  $(S, S,\delta_S)$  on $H$, 
there exists an element $(\mu_\gamma,V_\gamma)$ of $\Gamma$,  
with $\nu$  unitarily equivalent to $\mu_\gamma$ and 
$W$  unitarily equivalent to $V_\gamma$.
Let $\mu_S:=\oplus_{{}_\Gamma}\mu_\gamma$, 
$V_S:=\oplus{{}_\Gamma} V_\gamma$ and $H_S:=\oplus{{}_\Gamma} H_\gamma$.
The representation $\mu_S$ of $S$ is called the
{\it regular representation of $S$} and
 the pair $(\mu_S,V_S)$ is called the 
{\it regular covariant representation of $(S,S,\delta_S)$}.

Define $A_0(S)$ to be the pre-dual of the 
von Neumann algebra generated by $\mu_S(S)$ in $B(H_S)$.
For each $f\in A_0(S)$,
define $f_S\in S^*$
by  $f_S:=f\circ\mu_S.$
We may put operations on  $A_0(S)$ by
  pointwise addition and scalar multiplication, 
and multiplication  $\star$   
defined by 
\[f\star g(\mu_S(x)):=f_S\tensor g_S(\delta_S(x)),\]

\ni  \cite[Eq 1.1]{gl:qg}.
Since $\mu_S$ is part of a covariant pair,
$\ker\mu_S=\ker((\mu_S\tensor\mu_S)\circ\delta_S)$,
and  
  $f\star g$ is well-defined.
\ni To show that the product $f\star g$ is an element of $A_0(S)$,
we need to know that it is weakly continuous. 
Since $f\tensor g$ is weakly continuous,
it follows from the covariance of 
$(\mu_S,V_S)$ because
\[f\star g(\mu_S(x))
=f\tensor g(\id\tensor\mu_S(V_S))(\mu_S(x)\tensor 1)(\id\tensor\mu_S(V_S))^*).\]
The associativity    follows from  
the comultiplication identity  Equation (\ref{coact.identity}):
\begin{eqnarray*}
f\star (g\star h)(\mu_S(x))
&=&f\tensor (g\star h)(\mu_S\tensor\mu_S(\delta_S(x)))\\
&=& f_S\tensor g_S\tensor h_S(\id\tensor\delta_S(\delta_S(x)))\\
&=&f_S\tensor g_S\tensor h_S(\delta_S\tensor\id(\delta_S(x)))\\
&=&(f\star g)\tensor h(\mu_S\tensor\mu_S(\delta_S(x)))
=(f\star g)\star h(\mu_S(x)).
\end{eqnarray*}

\ni Hence $A_0(S)$ is an algebra.
Let 
\[M:=\{f\in A_0(S)|\id\tensor f_S(W)=0\ 
\forall\  {\rm corepresentations\  }W{\rm \ of\ }S\}.\]

\ni Using Equation (\ref{corep.id}), we see that $M$ is an ideal in $A_0(S)$:
\begin{eqnarray}\label{eq:multip}
\id\tensor (f\star g)_S(W)
&=&\id\tensor f_S\tensor g_S(\id\tensor\delta_S(W))\nonumber\\
&=&\id\tensor f_S\tensor g_S(W_{12}W_{13}) 
=\id\tensor f_S(W)\id\tensor g_S(W).
\end{eqnarray}
\ni
Let $A(S)$ be the quotient of  $A_0(S)$ by $M$. 

\m
A Hopf \c algebra $S$ is called {\it nondegenerate} if the ideal
$M$ is zero.

\begin{remark}
It is tempting to suggest that 
a Hopf \c algebra $S$ is nondegenerate if and only if $S$ is bi-simplifiable.
But we have only been able to show one direction of the implication,
  that
$M=\{0\}$ implies $S$ is bi-simplifiable.

Katayama showed that, for a coaction $\delta$ of a locally compact group $G$ on
a \c algebra $A$, $\delta$ is nondegenerate in Landstad's sense
if and only if $\delta(A)(1\tensor \CG)=A\tensor \CG$ \cite[Th 5]{kat}
(Landstad's nondegeneracy is that if
$(f\tensor g)\circ\delta=0$ for all $g\in B(G)$,
then $f=0$ in $A^*$).
For a comultiplication $\delta_S$
one can show that it is bi-simplifiable if and only if 
it is nondegenerate in Landstad's sense:
$(f_S\tensor g_S)\circ\delta_S=0$ for all $g\in A(S)$ 
implies $f=0$ in $A(S)$,
and,
$(g_S\tensor f_S)\circ\delta_S=0$ for all $g\in A(S)$ 
implies $f=0$ in $A(S)$.
We show that if $M=\{0\}$, then $\delta_S$ is nondegenerate in Landstad's 
sense.

Suppose $f\in A(S)$,
$M=\{0\}$ and $(f_S\tensor g_S)\circ\delta_S=0$ for all $g\in A(S)$.
Then    
$\id\tensor f_S\tensor g_S(\id\tensor\delta_S(W))=0$ for all corepresentations
$W$ of $S$. 
By Equation (\ref{corep.id}) this means
$\id\tensor f_S\tensor g_S(W_{12}W_{13})
=\id\tensor f_S(W)\id\tensor g_S(W)=0$ for all $g$ and $W$. 
Thus $\id\tensor f_S(W)=0$ for all $W$, and $f\in M$.
Since $M$ is zero, $f=0$ in $A(S)$.
\end{remark}


A Hopf \c algebra $S$ is 
{\it co-involutive} if
   for each $f\in A(S)$,
there exists a functional  $f^*\in A(S)$ such that 
\[\id\tensor (f^*)_S(W)=(\id\tensor f_S(W))^*,\] 
\ni for all corepresentations $W$ of $S$.
Such an $f^*$ is the unique functional in $A(S)$ satisfying this
relation because elements in the ideal $M$ are
zero in $A(S)$.


\m\ni{\bf Example 1(c).}
Consider   $A$ from Example 1(a).
Since the comultiplication is trivial,
its only corepresentation is 
$1\tensor 1$, and
 the pair $(\pi,1\tensor 1)$ is covariant 
for every representation $\pi\colon A\to B(H)$. 
So $\Gamma$   has an element from each equivalence class
of representations of $A$
and $\mu_S$ is the universal representation.
Thus $A_0(S)$ is simply the dual space $A^*$.
The ideal $M$ in this case is 
\[\{f\in A^*|\id\tensor f(1\tensor 1)=0\}
=\{f\in A^*|f(1)=0\}.\]
\ni We can realize $A(S)$ as the complex numbers ${\Bbb C}$ via the isomorphism
$f+M \mapsto f(1)$.
Complex conjugation satisfies the
requirements of a   co-involution, so
$A$ is a co-involutive Hopf \c algebra,
 but it is not nondegenerate.


\m\ni{\bf Example 2(c).}
Consider the Hopf \c algebra $C_0(G)$.
From Examples 2(a) and (b), we know that
$\alpha_G$ is a coaction of $C_0(G)$ on itself.
The covariant representations of 
$(C_0(G),C_0(G),\alpha_G)$ are the same as those of the dynamical system
$(C_0(G),G,\sigma)$, where $\sigma$ is   action by right translation.
By \cite{rie:compacts},
the covariant representations of this dynamical system are in 
one-to-one correspondence with the representations of $K(L^2(G))$.
Now, $K(L^2(G))$ is a simple, Type I \c algebra,
 so it has only one equivalence class of representations,
namely the faithful ones.
So, in this example, there is only one element in $\Gamma$.
The  covariant representation  $(M,\lambda)$ of $(C_0(G),G,\sigma)$ 
 has $M$
faithful,
and we choose it to be the element of $\Gamma$.
That is, $\mu_S:=M$ and $V_S:=\lambda\tensor\id(v_G)$.
Since $M$ is a faithful representation  
the weak closure of $M(C_0(G))$ is $L^\infty(G)$, 
whose predual  is $L^1(G)$.

 In Example 2(a) we showed
 that every corepresentation of $C_0(G)$
is of the form 
$U\tensor \id(v_G)$,
where $U$ is a
*-representation of $\CG$.
Thus the ideal $M$   is 
\[\{f\in L^1(G)|\id\tensor f(U\tensor \id(v_G))=0\ \forall\ U\}
=\{f\in L^1(G)|U(f)=0\ \forall\ U\}
=\{0\},\]
\ni  so   $A(C_0(G))$ is  $L^1(G)$.
Using Equation (\ref{eq:rock2}),  for $f\in L^1(G)$, we have
\begin{eqnarray*}
\id\tensor f^*(U\tensor \id(v_G))
&=&U(\id\tensor f^*(v_G))
=U(f^*)=U(f)^*\\
&=&[U(\id\tensor f(v_G))]^*
=[\id\tensor f(U\tensor \id(v_G))]^*.
\end{eqnarray*} 
\ni This shows that $C_0(G)$ is a nondegenerate co-involutive Hopf \c algebra.

 This next calculation shows that the operations defined 
on $A(C_0(G))=L^1(G)$
 give group multiplication between the
point masses in $L^1(G)$:
\[\varepsilon_s\star\varepsilon_t(f)
=\varepsilon_s\tensor\varepsilon_t(\alpha_G(f))
=\varepsilon_s\tensor\varepsilon_t((u,v)\mapsto f(uv))
=f(st)=\varepsilon_{st}(f).\]
\ni 
When an element of $L^1(G)$ is viewed as a measure on $G$
the multiplication defined here corresponds to
convolution of measures
\cite[7.1.2]{pedbook}.


\m\ni{\bf Example 3(c).}
Consider  the Hopf \c algebra $\CG$.
The argument here is similar to that in Example 2(c).
Again the covariant representations of
$(\CG,\CG,\delta_G)$ are 
precisely the covariant representations
of the cosystem 
$(\CG,G,\delta_G)$, which are
in one-to-one correspondence with the representations of 
$K(L^2(G))$ \cite{rie:compacts}, \cite[Ex  2.9]{iain:coact}.
So,   we can 
 choose
$(\mu_S, V_S)$ to be $(\lambda, M\tensor\id(w_G))$.
The group von Neumann algebra $VN(G)$
is defined to be the weak closure of the
image of $\CG$ under $\lambda$,
 that is
$\overline{\lambda(\CG)}{}^w=VN(G)$.
The predual of $VN(G)$ is the Fourier algebra $A(G)$ 
\cite[Th 3.10]{eymard} \cite[7.2.2]{pedbook}.

We showed in Example 3(a)
 that every corepresentation of $\CG$
is of the form 
$\mu\tensor \id(w_G)$,
where $\mu$ is a
*-representation of $C_0(G)$. 
Thus the ideal $M$  is 
\[\{f\in A(G)|\id\tensor f(\mu\tensor \id(w_G))=0\ \forall\ \mu\}
=\{f\in A(G)|\mu(f)=0\ \forall\ \mu\}
=\{0\},\]
\ni  so  in this example $A(\CG)$ is  $A(G)$ 
(hence the notation).

Using Equation (\ref{eq:rock}), we show that $\CG$ is a nondegenerate
co-involutive Hopf \c algebra:
%
\[\id\tensor f^*(\mu\tensor \id(w_G))
=\mu(\id\tensor f^*(w_G))
=\mu(f^*)=\mu(f)^*
=\id\tensor f(\mu\tensor \id(w_G))^*.\]

This next calculation shows that the multiplication defined  on $A(G)$
gives pointwise multiplication:
\[f\star g(s)
=f\tensor g(\delta_G(s))
= f\tensor g(i_G(s)\tensor i_G(s))
=f(s)g(s)=fg(s).\]

\begin{remark}
The weakness in the definition of the co-involutive property
is that one must know all the corepresentations of a 
Hopf \c algebra before 
you can check whether or not it is co-involutive.
\end{remark}

\m\ni{\bf Example 6(c).}
Let $V$ be an amenable multiplicative unitary which is part of a Kac triplet
\cite[Defn 6.4, p485]{bsk}.
We now show that Hopf \c algebra $S_V$ (see Example 6(a))
is co-involutive.
Firstly,  
$V$ is a corepresentation of $S_V$ and
  $(\id,V)$ is a covariant representation of
$(S_V,S_V,\delta)$, by a similar calculation as   in 
Example 6(b).
Baaj and Skandalis show that if $V$  is part of a Kac triplet,
the crossed product $S_V\times _\delta S_V$
is the compact operators on $H$ \cite[Prop 6.3, p485]{bsk},
and thus there is only one unitary equivalence class or
covariant representations of $(S_V,S_V,\delta)$.
Thus 
 $\mu_S=\id$ and $V_S=V$.
Define $\phi\colon A_0(S_V)\to \hat A(V)$
by
 $\phi(f)=\id\tensor f(V)$.
This is well-defined because if $f=0$,
this means that $f$ is in the pre-annihilator of $S_V$,
so that $f=0$ on $S_V$.
Thus 
for all $g\in B(H)_*$,
$f(g\tensor\id(V))=g(\id\tensor f(V))=0$,
which implies $\id\tensor f(V)=0$ in $\hat A(V)$.

The ideal $M=\ker\phi$ in $A_0(S_V)$ is zero by the following argument.
Suppose $f\in A_0(S_V)$ is positive and
$\id\tensor f(W)=0$ for all corepresentations $W$ of $S_V$.
In particular then $\id\tensor f(V)=0$.
So for all $g\in B(H)_*$,
$g(\id\tensor f(V))=f(g\tensor\id(V))=0$ and
  $f=0$ on $S_V$.
Now $S_V$ acts nondegenerately on $H$ \cite[Rem3.11(b)]{bsk},
so $S_V K(H)=K(H)$.
Since $f$ is positive,  $f=0$ on  $K(H)$,
and it is  weakly continuous, so it must be zero on the weak closure  $B(H)$.
Since $\phi$ is onto by the definition of $\hat A(V)$,
$\phi$ is both one-to-one and onto, and thus
  $A(S_V)$ is $\hat A(V)$.

Take $f\in A(S_V)$
and define $f^*$ to be the unique functional satisfying
$\id\tensor f^*(V)=(\id\tensor f(V))^*$.
This gives an element in 
$\hat A(S_V)$ because $V$ is regular so $\hat A(V)$ is a *-algebra 
\cite[Prop 3.5]{bsk}.

A multiplicative unitary is {\it amenable} if
$\hat S_V=\hat S_p$ \cite[Def 1.17]{ng:cc}, \cite[Prop 5.5]{blanch}.
So, since $V$ is an amenable multiplicative unitary,
it follows from \cite[Lemma 2.6]{ng:cc} that
every corepresentation $W$ of $S_V$ gives rise to a *-representation 
$\mu_{{}_W}$ of $\hat S_V$ satisfying
$W=\mu_{{}_W}\tensor\id(V)$.
Thus the following calculation shows that
  $S_V$ is co-involutive:
\begin{eqnarray*}
\id\tensor f^*(W) &=&\id\tensor f^*(\mu_{{}_W}\tensor\id(V))\\
&=&\mu_{{}_W}(\id\tensor f^*(V)) 
=\mu_{{}_W}((\id\tensor f(V))^*)\\
&=&[\mu_{{}_W}((\id\tensor f(V)))]^*
=[(\id\tensor f(W))]^*.
\end{eqnarray*}

\ni As an example,
the multiplicative unitary associated to  
  $Sl_\mu (2)$ defined by Woronowicz 
is amenable and part of a Kac triplet \cite[Ch 7]{blanch},
and so gives rise to a co-involutive Hopf \c algebra.

\begin{remark}
From the example we can see that the definition of 
co-involutive given here is suitable for
full Hopf \c algebras rather than reduced.
\end{remark}

\begin{lemma}\label{A(S).alg}
Let $S$ be  a co-involutive  Hopf \c algebra.
Then $A(S)$ is a *-algebra with 
  involution  satisfying
$\id\tensor f_S(W)^*=\id\tensor (f^*)_S(W)$,
 for all corepresentations $W$ of $S$. 
	%
\end{lemma}
\begin{pf}
We need to verify that $f^{**}=f$
 and  $(f\star g)^*=g^*\star f^*$:
\[\id\tensor (f^{**})_S(W)=\id\tensor (f^*)_S(W)^*
=(\id\tensor f_S(W))^{**}=\id\tensor f_S(W),\]
\ni so that  $(f^{**})_S=f_S$.
Also, 
\begin{eqnarray*}
\id\tensor ((f\star g)^*)_S(W)
&=&\id\tensor (f\star g)_S(W)^*
=[\id\tensor f_S(W)\id\tensor g_S(W)]^*\\ 
&=&\id\tensor g_S(W)^*\ \id\tensor f_S(W)^*\\
&=&\id\tensor (g^*)_S(W)\ \id\tensor (f^*)_S(W)
=\id\tensor (g^*\star f^*)_S(W),
\end{eqnarray*}
\ni so that $((f\star g)^*)_S=(g^*\star f^*)_S$.
We will be finished them both when we show that 
$f_S=g_S$ in $S^*$ implies $f=g$ in $A(S)$,
or equivalently,
$f_S=0$ in $S^*$ implies $f=0$ in $A(S)$. 
If $f_S=0$ in $S^*$,
then $\id\tensor f_S(W)=0$ for all corepresentations $W$ of $S$.
 This means $f$ is in the ideal $M$ of $A_0(S)$,
and thus is zero in the quotient $A(S)$.
\end{pf}


A Banach space $Z$ is a {\it nondegenerate  Banach $S$-module}
is a Banach $S$-module if $S$ is a \c algebra
and $Z$ is a Banach $S$-module \cite[32.14]{hr:II}
such that 
there exists an approximate identity $\{e_\lambda\}$ in $S$ satisfying
$e_\lambda\cdot z\to z $ for all $z\in Z$.


\begin{lemma}\label{prop:A(S)module}
Let $S$ be a Hopf \c algebra. Then
there is an action of $S$ on $A(S)$ satisfying
$( f \cdot x)_S(y)=f_S(xy)$,
 and $A(S)$ is a nondegenerate Banach $S$-module
with norm
$\|f\|=\sup\{|f(\mu_S(x))|\colon x\in S, \|x\|\leq 1\}$.
\end{lemma}
\begin{pf}
Certainly $A_0(S)$ is complete with respect to the norm 
it inherits as a subspace of   
$\overline{\mu_S(S)}{}^w$ in $B(H_S)$,
that is,
$\|f\|_\mu:=\sup\{|f(y)|\colon y\in \mu(S), \|y\|\leq 1\}.$
There is an isometric isomorphism between
the subspace of functionals that annihilate $\ker\mu_S$
and  $\left(S/(\ker\mu_S)\right)^*$.
This implies that $\|f\|$ is equal to 
 $\|f\|_\mu$.
So $A_0(S)$ is complete with respect to the norm  $\|\cdot \|$.

To see that $A(S)$ is complete we need to show that
the ideal \newline $M:=\{f\in A_0(S)|\id\tensor f_S(W)=0\ 
\forall\  {\rm corepresentations\  }W{\rm \ of\ }S\}$ is closed. 
It suffices to show that each
$M_W:=\{f\in A_0(S)|\id\tensor f_S(W)=0\}$.
Well, $M_W$ is the kernel of the map 
$\mu_{{}_W}\colon A_0(S) \to B(H)$
satisfying 
$\mu_{{}_W}(f)=\id\tensor f_S(W)$,
so if we show that  $\mu_{{}_W}$ is continuous,
we are done:
\[\|f\|=\|f_S\|=\|\id\tensor f_S\|
=\textstyle \sup_{\|x\|\leq 1}\displaystyle\|\id\tensor f_S(x)\|
\geq\|\id\tensor f_S(W)\|
=\|\mu_{{}_W}(f)\|.\]

Notice that the norm of $f\in A(S)$ is simply
the norm of $f_S$ in $S^*$.
Using this, together with
 the fact  $S^*$ is a nondegenerate  Banach $S$-module
\cite[p753]{lprs},
one can show that  $A(S)$ is also nondegenerate Banach $S$-module.
\end{pf}


%
%
%
%
%
%
%
%
%


\section{Crossed products of Hopf \c algebras.}\label{sec:crosprod}

Let  $\delta$ be a nondegenerate coaction of a  co-involutive Hopf \c algebra
$S$ on a \c algebra $A$.
The triple $(A, S, \delta)$ is called  a  {\it Hopf system}.
%
%
A \hbox{{\it full crossed product}} for a Hopf system  $(A, S, \delta)$ 
is a \c algebra $B$,
together with nondegenerate *-homomorphism 
$j_A\colon A\to M(B)$ and a unitary 
$u_S\in M(B\tensor  S)$ satisfying:

\begin{itemize}
\item[(a)] the pair $(j_A, u_S)$
is a covariant homomorphism of
 $(A, S, \delta)$ into $M(B)$,
\item[(b)] for every covariant representation $(\pi,W)$ of $(A, S, \delta)$,
there is a  representation  $\pi\times W$ of $B$ 
such that
$(\pi\times W)\circ j_A=\pi$ and $(\pi\times W)\tensor\id(u_S)= W$, and
\item[(c)] 
$B=\clsp\{j_A(a)\id\tensor f_S(u_S)\colon a\in A, f\in  A(S)\}$.
\end{itemize}

\ni  Before we show that such a crossed product exists,
we will show that the vector space in (c) is actually
a \c algebra.
The proof is modelled on  
\cite[Lem 2.10]{iain:coact}.

\begin{lemma}\label{B.is.a.c*.alg}
Let   $(A, S, \delta)$ be a  Hopf system, 
$C$ be a \c algebra,
and $(j_A,u_S)$ be a covariant homomorphism into $M(C)$.
Then 
 
\centerline{
$B:=\clsp\{j_A(a)\id\tensor f_S(u_S)\colon a\in A, f\in A(S)\}$}
\ni is a \c algebra.
\end{lemma}
\begin{pf} 
We begin by showing 
  that $B$ is  closed under multiplication.
For $f,g\in A(S)$,
\begin{eqnarray}\label{eq:phi.mult}
\id\tensor f_S(w_S)\id\tensor g_S(u_S)
&=&\id\tensor f_S\tensor g_S(u_{12}u_{13})\nonumber\\
&=&\id\tensor f_S\tensor g_S(\id\tensor \delta_S(u_S))\nonumber\\
&=&\id\tensor (f\star g)_S(u_S).
\end{eqnarray}

Take $a\in A$ and  $f\in A(S)$.
From the covariance of  $(j_A,u_S)$ we have that \newline
$\id\tensor f_S(u_S) j_A(a)
=\id\tensor f_S(j_A\tensor\id (\delta(a)) u_S).$
\ni Since $A(S)$ is a nondegenerate Banach $S$-module,
Lemma \ref{prop:A(S)module},
we can apply the Cohen factorization theorem \cite[p268]{hr:II} 
to see that
 $f\in A(S)$ can be expressed in the form $g\cdot x $,
where $g\in A(S)$, $x\in S$.
Thus
\begin{eqnarray*}
\id\tensor f_S(j_A\tensor\id (\delta(a)) u_S)
&=&\id\tensor (g\cdot x )_S(j_A\tensor\id (\delta(a)) u_S)\\
&=&\id\tensor g_S(j_A\tensor\id ( 1\tensor x.\delta(a)) u_S).
\end{eqnarray*}
\ni
Since the coaction satisfies $\delta(A)1\tensor S=A\tensor S$,
 we can approximate  
$(\delta(a))(1\tensor x)$
by a sum of the form $\sum a_i\tensor x_i$, so that
\begin{eqnarray*}
\id\tensor f_S(j_A\tensor\id (\delta(a)) u_S)
&\sim & \id\tensor g_S(j_A\tensor\id 
  (\textstyle\sum \displaystyle a_i\tensor x_i) u_S)\\
&=& \id\tensor g_S(\textstyle\sum  \displaystyle j_A(a_i)\tensor x_i. u_S)\\
&=& \textstyle\sum \displaystyle j_A(a_i)\, (\id\tensor (g\cdot x_i)_S  )(u_S).
\end{eqnarray*}
\ni Thus $B$ is closed under multiplication.
We next  show that 
\begin{equation}\label{eq:adjoint}
\id\tensor f_S(u_S)^*=\id\tensor f^*_S(u_S).
\end{equation}
\ni Let $\pi\colon B\to B(H)$ be a faithful representation.
Then $\pi\tensor\id(u_S)$ is a corepresentation of $S$ and
thus
\[[\pi(\id\tensor f_S(u_S))]^*
=[\id\tensor f_S(\pi\tensor\id(u_S))]^*
=\id\tensor f^*_S(\pi\tensor\id(u_S))
=\pi(\id\tensor f^*_S(u_S)).\]
\ni Equation (\ref{eq:adjoint}) follows from the fact that
 $\pi$ is a faithful *-homomorphism.
This equation shows that the adjoint of 
$j_A(a)\id\tensor f_S(u_S) $ 
 is in $B$. 
Hence $B$ is ${}^*$-closed.
\end{pf}

Let $(A, S,\delta)$ be a  Hopf system.
A covariant representation $(\pi, W)$ of $(A, S,\delta)$  
 on $H$ is 
{\it cyclic} if there exists a vector $\xi$ in $H$ 
such that
\[H=\clsp\{\pi(a)\id\tensor f_S(W)(\xi)\colon a\in A, f\in  A(S)\}\] 
\ni \cite[Defn 2.8(c)]{ng:cc}.
%
%
Let $(\mu_S,V_S)$  be the regular covariant representation of
$(S,S,\delta_S)$
on   $H_S$ (\S \ref{sec:dual.hopf}).
We know there are non-trivial
covariant representations of 
$(A,S,\delta)$ -- given a representation $\pi$ of $A$,
 $(\pi\tensor\mu_S\circ\delta,1\tensor V_S)$ is 
 covariant:
\begin{eqnarray*}
[((\pi\tensor \mu )\circ\delta)\tensor\id(\delta(a))] \, [1\tensor V]
&=&[(\pi\tensor \mu \tensor\id)
    \circ(\delta\tensor\id)(\delta(a))] \, [1\tensor V]\\
&=&[(\pi\tensor \mu \tensor\id)
  \circ(\id\tensor\delta_S)(\delta(a))] \, [1\tensor V]\\
&=&[\pi\tensor ((\mu \tensor\id)\circ\delta_S)(\delta(a))] \, [1\tensor V]\\
&=&[1\tensor V] \, [(\pi\tensor \mu (\delta(a)))\tensor 1].
\end{eqnarray*}
Define $\Ind\pi\colon A\times_\delta S\to B(H\tensor H_S)$
by $\Ind\pi:=(\pi\tensor\mu_S\circ\delta)\times(1\tensor V_S)$.
%
The following theorem is modelled on \cite[Prop  2.13]{iain:coact}.


\begin{theorem}\label{th:c.p.exists}
Let $(A,S,\delta)$ be  a Hopf system.
Then there is a 
full crossed product $(B,j_A,j_S)$
for $(A,S,\delta)$,
which is unique in the sense that if $(C,k_A,k_S)$ is another,
then there is an isomorphism $\phi$ of $B$ onto $C$ such that
$\phi\circ k_A=j_A$ and $\phi\circ k_S=j_S$.
The full crossed product is denoted by $A\times_\delta S$.
\end{theorem}
\begin{pf}
Let $\Gamma$ be a set of  cyclic covariant representations  
of   $(A, S,\delta)$
such that, for every   covariant representation $(\nu,W)$ of 
  $(A, S,\delta)$, 
there exists a member  $(\mu_\gamma,V_\gamma)$ of $\Gamma$,
with $\nu$  unitarily equivalent to $\mu_\gamma$ and 
$W$  unitarily equivalent to $V_\gamma$.
Let $j_S:=\oplus_{{}_\Gamma}\mu_\gamma$, 
$u_S:=\oplus_{{}_\Gamma} V_\gamma$ and $H:=\oplus_{{}_\Gamma} H_\gamma$.
Let $B$ be the closed linear span of 
$\{j_A(a)\id\tensor f_S(u_S)\colon a\in A, f\in  A(S)\}$.
By Lemma \ref{B.is.a.c*.alg}, $B$ is a \c algebra
and $j_A(a)$ and $\id\tensor f_S(u_S)$ are both multipliers of $B$.
As in \cite[p635]{iain:coact},
$(j_S,u_S)$ is a covariant homomorphism into $B$, so condition (a) is satisfied. Condition (c) holds by the definition of $B$.
  
As in \cite[Cor 2.12]{iain:coact},
it follows from Lemma \ref{B.is.a.c*.alg}
 that any covariant representation is
a direct sum of cyclic representations.
It is enough to check  condition (b) for a cyclic 
covariant representation $(\mu,V)$.
But  $(\mu,V)$ is equivalent to a member of $\Gamma$
and we can construct 
 $\mu\times V$ by compressing $B$ to the appropriate summand of $H$
(as in proof of \cite[Prop 2.13]{iain:coact}).
The uniqueness of the full crossed product follows 
immediately from its universal properties.
\end{pf}


\m\ni{\bf Example 4(c).}
Now we show that given a dynamical system
 $(A, G, \alpha)$,
the crossed product $A\times_\alpha G$ is a crossed product for the  Hopf system 
 $(A, C_0(G), \alpha)$.
The embeddings are $(k_A, k_G\tensor \id(v_G))$,
so condition (a) is satisfied.
Let $(\pi,W)$ be a covariant representation of
$(A, C_0(G), \alpha)$.
In Example 2(a)
we showed that there exists a  representation $V$ of $\CG$
satisfying
$V(s)=\id\tensor\varepsilon_s(W)$.

Since $W=V\tensor\id(v_G)$,
 $(\pi,V)$ is a covariant representation
of $(A,G,\alpha)$.
Thus $\pi\times V$ exists,
$(\pi\times V)\circ k_A=\pi$
and 
\[(\pi\times V)\tensor\id(k_G\tensor\id(v_G))
= V\tensor\id(v_G)=W,\]
\ni so condition (b) is satisfied.
It remains to verify condition (c).
Since $A(C_0(G))=L^1(G)$ is dense in $\CG$, we have 
\begin{eqnarray*}
A\times_\alpha C_0(G)
&=&\clsp\{k_A(a)\id\tensor z (k_G\tensor\id(v_G))
        \colon a\in A, z\in L^1(G)\}\\
&=& \clsp\{k_A(a) k_G(z) \colon a\in A, z\in L^1(G)\}
=A\times_\alpha G.
\end{eqnarray*}


\m\ni{\bf Example 5(c).}
For  a  cosystem $(A, G, \delta)$, 
the argument that $A\times_\delta G$ is a crossed product for the  Hopf system 
$(A,\CG, \delta)$ follows closely that of Example 4(c),
so that
%
%
%
%
%
%
%
$
A\times_\delta \CG
=\clsp\{j_A(a)j_{C_0(G)}(f)
        \colon a\in A, f\in A(G)\}$.

\begin{theorem}\label{th:cov.hom}
Let $(A,S,\delta)$ be a Hopf system, 
 $B$ be a \c algebra 
and $(\pi,W)$ be a covariant homomorphism into $M(B)$.
Then there exists a unique nondegenerate homomorphism
$\pi\times W\colon A\times_\delta S\to M(B)$
such that 
\[(\pi\times W)\circ j_A=\pi  \quad  and\quad 
 (\pi\times W)\tensor\id(u_S)=W.\]
\end{theorem}
\begin{pf} Represent $B$ faithfully on Hilbert
space via its universal
 representation $\phi\colon B\to B(H)=M(K(H))$,
so that
$(\bar\phi\circ\pi,\bar\phi\tensor\id(W))$ 
is a covariant representation of $(A,S,\delta)$.
Then 
there is a  representation
$\sigma:=(\bar\phi\circ\pi)\times(\bar\phi\tensor\id(W))$
of $A\times_\delta S$ on $H$ such that
$\sigma\circ j_A=\bar\phi\circ\pi  \and
\sigma\tensor\id(u_S)=\bar\phi\tensor\id(W).$
\ni Since the image of $\sigma$ is contained in the image of $\phi$,
$\bar\phi^{-1}\circ\sigma\colon A\times_\delta S\to M(B)$
is a homomorphism such that
$(\bar\phi^{-1}\circ\sigma)\circ j_A
=\bar\phi^{-1}\circ\bar\phi\circ\pi=\pi,
\ \and\ 
(\bar\phi^{-1}\circ\sigma)\tensor\id(u_S))=W.$
\ni The argument of    \cite[Lem 1.3]{PRII}
gives the nondegeneracy.
\end{pf}


\begin{remark}
Let $(A,S,\delta)$ be a Hopf system.
One would like to define the 
 reduced crossed product $A\times_{\delta,r}S$ 
as the  quotient of
$A\times_\delta S$ by $\ker(\Ind\pi)$
where $\pi$ is a faithful representation of $A$.
But for this to work, we need to know that $\Ind$ 
is well-defined on ideals.
We have not been able to show this is the case without
further assumptions on $S$.
This will be discussed further in \cite{nil:hopfresind}.
\end{remark}



\section{The Dual Hopf \c algebra.}\label{subsec:can.unit}

We  noted in Example \ref{ex:hopf.C^*-alg=hopf.system}
that  any Hopf \c algebra $S$ coacts trivially
on the complex numbers $\Bbb{C}$.
For a co-involutive Hopf \c algebra $S$,
define
\[\hat S:={\Bbb C}\times_{\id} S.\]
From the universal properties of the crossed product
 we have canonical unitaries 
$w_S\in M(\hat S\tensor S)$  and  
  $v_S:=\Sigma(w_S)\in M(S\tensor \hat S)$, 
\ni  satisfying
\begin{equation}\label{eq:wS.condit}
\id\tensor\delta_S(w_S)=(w_S)_{12}(w_S)_{13}
         \in M(\hat S\tensor S\tensor S),\ \and
\end{equation}
\[\delta_S\tensor\id(v_S)=(v_S)_{13}(v_S)_{23}
            \in   M(S\tensor S\tensor \hat S).\]

\begin{theorem}\label{hat.S.is.a.hopf.alg}
Let  $S$ be a  co-involutive Hopf \c algebra.
Define 
\[\psi\colon A(S)\to  \hat S \ \ by\ \ 
\psi(f):=\id\tensor f_S(w_S),\]
\ni where $w_S$ is the canonical unitary in $M(\hat S\tensor S)$.
Then $\psi$ is an   injective *-homomor\-phism
from $A(S)$ into $\hat S$ with dense range. 

If $W$ is a corepresentation  of $S$,
 then $\nu_{{}_W}$,  
 defined by $\nu_{{}_W}(\psi(f)):=\id\tensor f_S(W)$,
  is a  *-representation of $\hat S$
such that $\nu_{{}_W}\tensor\id(w_S)=W$. 

Furthermore, $\hat S$ is a Hopf \c algebra with comultiplication 
$\delta_{\hat S}\colon \hat S\to M(\hat S\tensor \hat S)$ satisfying
\vspace*{-0.1cm}
\begin{equation}\label{eq:hat.coact}
\delta_{\hat S}(\psi(f))
   =\id\tensor  \id\tensor f_S((w_S)_{13}(w_S)_{23}).
\end{equation}
\end{theorem}
%
%
\begin{pf} The operations on $A(S)$ are defined in \S\ref{sec:dual.hopf}.
Equation (\ref{eq:phi.mult}) shows that  $\psi$   is  multiplicative.
 The norm on   
$\hat S={\Bbb C}\times_{\id} S$
is the supremum of the covariant representations of 
$({\Bbb C},S,\id)$,
which are exactly the  corepresentations of $S$, so  
\begin{eqnarray*}
\|\psi(f)\|&=&\|\id\tensor f_S(w_S)\|\\
&=&\sup_W\{\|\id\times W(\id\tensor f_S(w_S))\|
        \colon  W \mathrm{\ is\ a\ corepresentation\ of\ }S\}\\
&=&\sup_W\{\|\id\tensor f_S((\id\times W)\tensor\id(w_S))\|\}
=\sup_W\{\|\id\tensor f_S(W)\|\}.
\end{eqnarray*}
\ni Now if $\psi(f)=0$, the above shows that
$\id\tensor f_S(W)=0$ for all corepresentations $W$ of $S$,
which means  that $f=0$ in $A(S)$.
Thus $\psi$ is  injective.
The image of $\psi$ is dense in $\hat S$ because
$\hat S$ is spanned by 
 $\{\id\tensor f_S(w_S)\colon f\in A(S)\}$.

Let $W$ be a corepresentation  of $S$ on $H$;
 define $\nu_{{}_W}\colon \hat S \to B(H)$ 
by 
\[\nu_{{}_W}(\psi(f)):=\id\tensor f_S(W).\]
\ni A calculation like that of Equation (\ref{eq:phi.mult})
shows that it is multiplicative,
and it is involutive because $S$ is a co-involutive Hopf \c algebra.
It is norm decreasing since the norm in $\hat S$
is the supremum of the corepresentations of $S$.
Since the slice maps $\id\tensor f$ for $f\in S^*$ separate
points, it follows from the definition of $\nu_{{}_W}$ that
  $\nu_{{}_W}\tensor \id(w_S)=W$.

To show that there exists a nondegenerate *-homomorphism
 $\delta_{\hat S}$ satisfying Equation (\ref{eq:hat.coact}),
we will apply Theorem \ref{th:cov.hom}.
It will suffice to show that there exists
a unitary $u$ in 
$M(\hat S\tensor\hat S\tensor S)$
such that
$\id\tensor\id\tensor\delta_S(u)=u_{123}u_{124}$.
A routine calculation shows that 
$u:=w_{13}w_{23}$ does the job.
%
%
%
It follows from the definition of 
$\delta_{\hat S}$ that 
\begin{equation}\label{eq:delta.S.hat}
\delta_{\hat S}\tensor\id(w_S)=(w_S)_{13}(w_S)_{23}
\in M(\hat S\tensor \hat S \tensor S):
\end{equation}
The following shows that the comultiplication identity is satisfied:
\begin{eqnarray*}
\id\tensor\delta_{\hat S}(\delta_{\hat S}(f))
&=& \id\tensor\delta_{\hat S}(\id\tensor  \id\tensor f(w_{13}w_{23}))\\
&=& \id\tensor\id\tensor f
       (\id\tensor \delta_{\hat S}\tensor\id(w_{13}w_{23}))\\\
&=& \id\tensor\id\tensor\id\tensor f
    (w_{14} w_{24}w_{34})\\
&=& \id\tensor\id\tensor\id\tensor f
    (\delta_{\hat S}\tensor\id\tensor\id(w_{13} w_{23}))\\
&=& \delta_{\hat S}\tensor\id 
       (\id\tensor\id\tensor f (w_{13} w_{23})) 
= \delta_{\hat S}\tensor\id(\delta_{\hat S}(f)).  
\end{eqnarray*} 

\vspace*{-.5cm}

\end{pf}


\begin{remark}
Note that if $V$ and $W$ are unitarily equivalent corepresentations
of $S$, then the representations $\nu_{{}_V}$ and 
$\nu_{{}_W}$ of $\hat S$ are  unitarily equivalent.
\end{remark}


 From now on  the injection
$\psi\colon  A(S) \to \hat S$ given by Theorem \ref{hat.S.is.a.hopf.alg}
will  be used implicitly.
So, for example,
we can consider $f\in A(S)$ to be an element of $\hat S$,
and the canonical unitary $w_S$  
  in $M(\hat S\tensor S)$ satisfies
\begin{equation}\label{eq:can.unitary}
 \id\tensor f_S(w_S)=f\in A(S).
\end{equation}
\ni For $f\in S^*$,
$ \id\tensor f_S(w_S)=f\in M(\hat S)$.

The Hopf \c algebra $\hat S$ is called the
{\it dual Hopf \c algebra for $S$}.
It follows from Equation (\ref{eq:delta.S.hat}) that
\begin{equation}\label{eq: v_s.corep}
\id \tensor\delta_{\hat S}(v_S)
= (v_S)_{12} (v_S)_{13}
\in M(S\tensor \hat S\tensor \hat S).
\end{equation}
\ni Thus given a representation $\nu$ of $S$,
$\nu\tensor\id(v_S)$ is a corepresentation of $\hat S$.

\m

 Let $(A, S,\delta)$ 
be a  Hopf system.
By definition there is a 
(not necessarily faithful) *-homomorphism 
$j_A\colon A\to M(A\times_\delta S)$.
Now define 
$j_{\hat S}\colon \hat S\to M(A\times_\delta S) 
\ \ \by \ \ j_{\hat S}(f):=\id\tensor f(u_S),$
\ni where $u_S$ is the canonical unitary in 
$M((A\times_\delta S)\tensor S)$.
Calculations just like 
those in Theorem \ref{hat.S.is.a.hopf.alg}
show  that $j_{\hat S}$ is an injective  *-homomorphism.
Thus  
\[A\times_\delta S=\clsp\{j_A(a)j_{\hat S}(f)\colon a\in A, f\in A(S)\}.\]

\medskip\ni{\bf Example 1(b).}
In the trivial example $A$, $\hat A={\Bbb C}$,
so the trivial crossed product by $A$ of any \c algebra $B$
is nothing more than  $B$.


\medskip\ni{\bf Example 2(d).}
It was shown in Example 2(a)
  that  
$C_0(G)$ is a   Hopf \c algebra with comultiplication $\alpha_G$.  
It is well known that the dynamical system crossed product  
${\Bbb C}\times_{\id} G$ is isomorphic to $\CG$,
 so by Theorem \ref{hat.S.is.a.hopf.alg} 
and Example 4(c),
  the dual Hopf \c algebra 
$\widehat{C_0(G)}$ is $\CG$.
The comultiplication  
 specified
  is in fact  $\delta_G$:
\begin{eqnarray*}
\delta_{\hat S}(\varepsilon_s)
&=&\id\tensor  \id\tensor \varepsilon_s((v_G)_{13}(v_G)_{23})\\
&=&\id\tensor  \id\tensor \varepsilon_s 
            ((u\mapsto i_G(u))(v\mapsto i_G(v)))\\
&=&i_G(s)\tensor i_G(s)=\delta_G(s).
\end{eqnarray*}
 

\medskip\ni{\bf Example 3(d).}
It was shown in  Example 3(a)
  that  the full group \c algebra $\CG$ is a   Hopf \c algebra
 with comultiplication $\delta_G$. 
Since
${\Bbb C}\times_{\id}  G\cong C_0(G)$,
  the dual Hopf \c algebra $\widehat\CG$ is $C_0(G)$
with comultiplication  
   $\alpha_G$:
\begin{eqnarray*}
\delta_{\hat S}(f)
&=&\id\tensor  \id\tensor f((w_G)_{13}(w_G)_{23})\\
&=&\id\tensor  \id\tensor f    ((s\mapsto i_G(s))(t\mapsto i_G(t)))\\
&=&\id\tensor  \id\tensor f ((s,t)\mapsto   i_G(s)i_G(t))\\
&=&\id\tensor  \id\tensor f ((s,t)\mapsto  i_G(st))
=(s,t)\mapsto   f(st)=\alpha_G(f).
\end{eqnarray*}


\medskip\ni{\bf Example 6(d).}
Let $V$ be an amenable multiplicative unitary which is part of a Kac triplet.
We argued  in   Example 6(c)  that $\hat S_V$ is the \c algebra whose
representations are in one-to-one correspondence with the
corepresentations of $S_V$.
The formula  given for the comultiplication in 
Theorem \ref{hat.S.is.a.hopf.alg} agrees with that given by
Baaj and Skandalis in \cite[cor A.6]{bsk}, and
thus $\hat S_V$ is the same as the dual   constructed here.


\section{The double dual.}\label{sec:doubledual}

Let $S$ be a co-involutive Hopf \c algebra
and let  $(\mu,V)$ be a 
    covariant representation of $(S, S, \delta_S)$ on $H$.
Define   
\begin{eqnarray*}
\hat\mu\colon \hat S\to B(H_S)&\ \by\ &\hat\mu(f):=\id\tensor f(V), \and \\
\hat V\in M(K(H_S)\tensor \hat S)&\ \by\ &\hat V:=  \mu\tensor\id(v_S).
\end{eqnarray*}
\ni  By Theorem  \ref{hat.S.is.a.hopf.alg},
  $\hat\mu$ is a *-representation of $\hat S$.
Using 
Equation (\ref{eq: v_s.corep}) we can show that 
 $\hat V$ is a corepresentation of $\hat S$:
\begin{eqnarray*}\id \tensor\delta_{\hat S}(\hat V)
&=&\id \tensor\delta_{\hat S}(\mu\tensor\id(v_S))\\
&=&\mu\tensor\id(\id \tensor\delta_{\hat S}(v_S))
=\mu\tensor\id((v_S)_{12} (v_S)_{13})\\
&=&\mu\tensor\id((v_S)_{12})\mu_S\tensor\id( (v_S)_{13})
=(\hat V)_{12}(\hat V)_{13}.
\end{eqnarray*}


Now, it is natural to ask whether or not
$(\hat\mu, \hat V)$ 
is a covariant representation of   
$(\hat S, \hat S, \delta_{\hat S})$.
Both Quigg and Raeburn have sought such results
\cite[Ex 2.9]{iain:coact}, \cite[Prop 2.5]{qui:fcpd}.
But it doesn't always work. 
For example,
it doesn't work for 
$(\CG,\CG, \delta_G)$ and $(\lambda, M\tensor\id(w_G))$,
since   $(M, \lambda\tensor\id(v_G))$
is not covariant for $(C_0(G),C_0(G),\alpha)$.
In order to construct the dual of a Hopf \c algebra 
we need to know that there exists at least one covariant representation,
and as we have just seen, $\hat S$ does
not automatically have them.

So,  
in order to build the double dual $\widehat{\hat S}$
we need to know a number of things about  $(\hat S, \hat S, \delta_{\hat S})$,
in particular, that it is a co-involutive Hopf \c algebra.
For now, we suppose that $S$ is a Hopf \c algebra
such that $\hat S$ is co-involutive, 
and  denote its regular   covariant representation
by $(\mu_{\hat S}, V_{\hat S})$.
Is that enough to ensure that
$\widehat{\hat S}$ is isomorphic to $S$?


Given $x\in S$, 
evaluation at $x$ is a functional on $S^*$; 
  denote it by $\varepsilon_x$. 
This functional can be viewed as a functional on $A(S)$ 
if and only if $S$ if nondegenerate, that is,
 if the ideal $M$ as $A_0(S)$, defined in \S \ref{sec:dual.hopf},
is zero. 
In this case   $\varepsilon_x$ extends to a functional
on $\hat S$ and 
$f(x)=\varepsilon_x(f) 
=\varepsilon_x( \id\tensor f(w_S)) 
= f(\varepsilon_x\tensor \id(w_S)),$
\ni for  all  $f\in S^*$, (Equation (\ref{eq:can.unitary}))
which implies that $\varepsilon_x\tensor \id(w_S)=x$,
or equivalently,
\begin{equation}\label{eq:can.unitary2}
 \id\tensor \varepsilon_x(v_S)=x.
\end{equation}

One way to show that the double dual   
is isomorphic to $S$, would be to show that 
 $S$ is a crossed product for the Hopf system
$({\Bbb C}, \hat S, \id)$.%
We would need to check the conditions of the definition of a crossed product.
For a homomorphism of ${\Bbb C}$ into $M(S)$ we just map
$z$ to $z1$. 
For the  unitary in  $M(S\tensor \hat S)$ we just use $v_S$.
This pair is certainly covariant, so condition (a) is satisfied.

Let $W$ be a corepresentation of $\hat S$ on $H$.
We would need to show that there exists a representation $\nu$ of $S$
such that $\nu\tensor\id(v_S)=W$.
Well, 
define $\nu_{{}_W}\colon S\to B(H)$ by
$\nu_{{}_W}(x):=\id\tensor \varepsilon_x(W)$,
where $\varepsilon_x$ is evaluation at $x$.
Then, if $S$ is nondegenerate, then 
$\id\tensor \varepsilon_x(\nu_{{}_W}\tensor\id(v_S))
     =\nu_{{}_W}(\id\tensor \varepsilon_x(v_S))
     =\nu_{{}_W}(x)= \id\tensor \varepsilon_x(W),$
\ni for all $x\in S$  (Equation (\ref{eq:can.unitary2})).
Since $A(S)$ is dense in $\hat S$,
the point evaluation functionals $\varepsilon_x$
 are sufficient to separate points of $\hat S$.
Thus functionals of the form $\id\tensor\varepsilon_x$ separate points
of $M(K(H)\tensor \hat S)$,
   $\nu_{{}_W}\tensor\id(v_S)=W$ and condition (b) is satisfied.

For condition (c) we need to show that 
$S$ is equal to   
$\clsp\{\id\tensor g_{\hat S}(v_S)\colon g\in A(\hat S)\}$.
This is the problem.
We have nothing that relates the covariant pairs
of $(S,S,\delta_S)$ 
to the covariant pairs
of $(\hat S,\hat S,\delta_{\hat S})$,
so we have no way of relating the elements of $A(\hat S)$ to elements of $S$.

To be able to show that
the double dual $\widehat{\hat S}$  is isomorphic to $S$,
 we need a 
prescribed correspondence between the
 covariant pairs
of $(S,S,\delta_S)$ 
and the covariant pairs
of $(\hat S,\hat S,\delta_{\hat S})$.
For this we will define a Kac system,
which will be  modelled on the relationship between the representations, 
 $\lambda$ and $M$,
associated
to a locally compact group $G$,
and the unitary operator 
$S\in B(L^2(G))$   defined  by  $S(\xi)(s)=\xi(s^{-1})$.
There are some similarities between our approach and that of 
Baaj and Skandalis \cite[\S 6]{bsk}.
This issue will be investigated in \cite{nil:hopfresind}.

\providecommand{\bysame}{\leavevmode\hbox to3em{\hrulefill}\thinspace}


\begin{thebibliography}{10}

\bibitem{alnr}
S.~Adji, M.~Laca, M.~Nilsen, and I.~Raeburn, \emph{Crossed products by
  semigroups of endomorphisms and the {T}oeplitz algebras of ordered groups},
  Proc. Amer. Math. Soc. \textbf{122} (1994), 1133--1141.

\bibitem{baaj}
S.~Baaj, \emph{Representation reguliere du groupe quantique des deplacements de
  {W}oronowicz}, Asterisque (1995), no.~232, 11--48.

\bibitem{bsk}
S.~Baaj and G.~Skandalis, \emph{Unitaires multiplicatifs et dualit\'e pour les
  produits crois\'es de \c alg\`ebres}, Ann. scient. \'Ec. Norm. Sup.
  \textbf{26} (1993), 425--488.

\bibitem{blanch}
E.~Blanchard, \emph{D\'eformation de \c alg\'ebres de {H}opf}, Bull. Soc. Math.
  France \textbf{124} (1996), 141--215.

\bibitem{dixbook}
J.~Dixmier, \emph{\c algebras}, North-Holland, New York, 1977.

\bibitem{eymard}
P.~Eymard, \emph{L'alg\`ebre de {F}ourier d'un groupe localement compact},
  Bull. Soc. Math. France \textbf{92} (1964), 181--236.

\bibitem{gl:qg}
E.~C. Gootman and A.~J. Lazar, \emph{Quantum groups and duality}, Rev. Math.
  Phys. \textbf{5} (1993), 417--451.

\bibitem{hr:II}
E.~Hewitt and K.~A. Ross, \emph{Abstract harmonic analysis}, vol.~II,
  Springer-Verlag, Berlin, 1970.

\bibitem{it}
S.~Imai and H.~Takai, \emph{On a duality for \c crossed products by a locally
  compact group}, J. Math. Soc. Japan \textbf{30} (1978), 495--504.

\bibitem{kat}
Y.~Katayama, \emph{Takesaki's duality for a non-degenerate co-action}, Math.
  Scand. \textbf{55} (1985), 141--151.

\bibitem{lprs}
M.~B. Landstad, J.~Phillips, I.~Raeburn, and C.~E. Sutherland,
  \emph{Representations of crossed products by coactions and principal
  bundles}, Trans. Amer. Math. Soc. \textbf{299} (1987), 747--784.

\bibitem{man}
K.~Mansfield, \emph{Induced representations of crossed products by coactions},
  J. Funct. Anal. \textbf{97} (1991), 112--161.

\bibitem{mur}
G.~J. Murphy, \emph{\c algebras and operator theory}, Academic Press, New York,
  1990.

\bibitem{nt}
Y.~Nakagami and M.~Takesaki, \emph{Duality for crossed products of von
  {N}eumann algebras}, Springer-Verlag, Berlin, 1979.

\bibitem{ng:cc}
C.~K. Ng, \emph{Coactions and crossed products of {H}opf \c algebras}, Proc.
  Lond. Math. Soc. \textbf{72} (1996), 638--656.

\bibitem{nil:dual}
M.~Nilsen, \emph{Duality for full crossed products of \c algebras by
  non-amenable groups}, to appear, Proc. Amer. Math. Soc.

\bibitem{nil:coact}
\bysame, \emph{Full crossed products by coactions, {C}${}_0${(X)}-algebras and
  \c bundles}, to appear, J. London Math. Soc.

\bibitem{nil:hopfresind}
\bysame, \emph{Kac-{H}opf \c systems and the induction and restriction of
  ideals}, in preparation.

\bibitem{PRI}
J.~Packer and I.~Raeburn, \emph{Twisted crossed products of \c algebras}, Math.
  Proc. Camb. Phil. Soc. \textbf{106} (1989), 293--311.

\bibitem{PRII}
\bysame, \emph{Twisted crossed products of \c algebras {II}}, Math. Ann.
  \textbf{287} (1990), 595--612.

\bibitem{pedbook}
G.~K. Pedersen, \emph{\c algebras and their automorphism groups}, Academic
  Press, London, 1979.

\bibitem{qui:fcpd}
J.~C. Quigg, \emph{Full \c crossed product duality}, J. Austral. Math. Soc.
  (Series A) \textbf{50} (1991), 34--52.

\bibitem{qui:frc}
\bysame, \emph{Full and reduced \c coactions}, Math. Proc. Camb. Philos. Soc.
  \textbf{116} (1994), 435--450.

\bibitem{qr:ind}
J.~C. Quigg and I.~Raeburn, \emph{Induced \c algebras and {L}andstad duality
  for twisted coactions}, Trans. Amer. Math. Soc. \textbf{347} (1995),
  2885--2915.

\bibitem{iain:tak}
I.~Raeburn, \emph{On crossed products and {T}akai duality}, Proc. Edin. Math.
  Soc. \textbf{31} (1988), 321--330.

\bibitem{iain:coact}
\bysame, \emph{On crossed products by coactions and their representation
  theory}, Proc. London Math. Soc. \textbf{64} (1992), 625--652.

\bibitem{rie:compacts}
M.~Rieffel, \emph{On the uniqueness of the {H}eisenberg commutation relations},
  Duke Math. J. \textbf{39} (1972), 745--752.

\bibitem{rudin:gp}
W.~Rudin, \emph{Fourier analysis on groups}, Interscience Publishers, New York,
  1962.

\bibitem{tom:ten}
J.~Tomiyama, \emph{Applications of {F}ubini type theorem to the tensor product
  of \c algebras}, Tohoku Math. J. \textbf{19} (1967), 213--226.

\bibitem{marty:w*}
M.~E. Walter, \emph{${W}^{\displaystyle *}$\!-algebras and nonabelian harmonic
  analysis}, J. Funct. Anal. \textbf{11} (1972), 17--38.

\bibitem{wor:cmp}
S.~L. Woronowicz, \emph{Compact matrix pseudogroups}, Commun. Math. Phys.
  \textbf{111} (1987), 613--665.

\bibitem{wor:twi}
\bysame, \emph{Twisted {${SU(2)}$} group. {A}n example of a non-commutative
  differential calculus}, Publ. RIMS, Kyoto Univ. \textbf{23} (1987), 117--181.

\bibitem{yam:groupoid}
T.~Yamanouchi, \emph{Duality for generalized {K}ac algebras and
  characterization of finite groupoid algebras}, J. Algebra \textbf{163}
  (1994), 9--50.

\end{thebibliography}
\end{document}